\documentclass{article}

\usepackage{times}

\usepackage{amsmath}
\usepackage{amssymb}
\usepackage{amsfonts}
\usepackage{amsthm}
\usepackage{xfrac}
\usepackage{enumitem}
\usepackage{latexsym}
\usepackage{multirow}
\usepackage{pbox}

\usepackage{graphicx} 
\usepackage{subfigure} 

\usepackage{natbib}

\usepackage{algorithm}
\usepackage{algorithmic}

\usepackage{hyperref}

\theoremstyle{plain}
\newtheorem{thm}{Theorem}
\newtheorem*{thm*}{Theorem}
\newtheorem{lem}{Lemma}
\newtheorem{prop}{Proposition}
\newtheorem*{prop*}{Proposition}
\newtheorem{cor}{Corollary}

\theoremstyle{definition}
\newtheorem{defn}{Definition}

\theoremstyle{remark}

\newtheorem*{prf}{Proof}

\usepackage[accepted]{icml2014} 

\icmltitlerunning{Pseudocontractive Updates}

\begin{document} 

\twocolumn[
\icmltitle{Properties of Pseudocontractive Updates \break in Convex Optimization}

\icmlauthor{Patrick W. Gallagher}{patrick.w.gallagher@gmail.com}
\icmladdress{UC San Diego,
	    9500 Gilman Dr., La Jolla, CA 92093}
\icmlauthor{Zhuowen Tu}{ztu@ucsd.edu}
\icmladdress{UC San Diego,
            9500 Gilman Dr., La Jolla, CA 92093}

\icmlkeywords{pseudocontractive, firmly nonexpansive, machine learning, convex optimization, operator theory, averaged, inverse strongly monotone}

\vskip 0.3in
]

\begin{abstract} 
Many convex optimization methods are conceived of and analyzed in a largely separate fashion.  In contrast to this traditional separation,
this manuscript points out and demonstrates the utility of an important but largely unremarked common thread running through many 
prominent optimization methods. In particular, we show that methods such as successive orthogonal projection, gradient descent, projected gradient descent, the proximal-point method, forward-backward splitting, the alternating direction method of multipliers, and under- or over-relaxed variants of the preceding all involve updates that are of a common type --- namely, the updates satisfy a property known as pseudocontractivity. Moreover, since the property of  pseudocontractivity is preserved under both composition and convex combination, updates constructed via these operations from pseudocontractive updates are themselves  pseudocontractive. Having demonstrated that pseudocontractive updates are to be found in many optimization methods, we then provide a unified basic analysis of methods with pseudocontractive updates. Specifically, we prove a novel bound satisfied by the norm of the difference in iterates of pseudocontractive updates and we then use this bound to establish  that the error criterion $\left\Vert x^{N}-Tx^{N}\right\Vert ^{2}$ is $o(1/N)$ for any method involving pseudocontractive updates (where $N$ is the number of iterations and $T$ is the iteration operator).
\end{abstract}

\section{Introduction}
\label{introduction}

Optimization methods are traditionally viewed and analyzed largely separately from one another. In this manuscript,
we highlight a common feature, namely \emph{pseudocontractivity}, of the iterative updates of many methods, and we show how recognizing this 
commonality can lead to a unified analysis and deeper understanding of a variety of prominent convex optimization methods. 

We summarize our contributions as follows: (1) we highlight the essentially unremarked prevalence of pseudocontractive updates in many convex optimization methods; 
(2) we prove a novel result, Lemma \ref{lem:monotone-decreasing}, establishing a bound on the decrease of the norm of the displacement operator associated with any 
pseudocontractive operator;
(3) we use the Lemma \ref{lem:monotone-decreasing} result on monotone norm decrease to prove a novel result, Theorem \ref{thm:convergence-rate-stronger},  establishing that the error criterion $\left\Vert x^{N}-Tx^{N}\right\Vert ^{2}$ is $o(1/N)$ for any method involving pseudocontractive updates 
(where $N$ is the number of iterations and $T$ is the iteration operator).

In order to establish the results above, we also bring together a variety of previously widely dispersed results; the combination of these results turns out to greatly
enhance their collective utility and applicability.

\subsection{Related Work}
We believe that pseudocontractive operators provide a new and widely useful perspective in the context of machine learning (and largely unmentioned in optimization). Most standard optimization texts (e.g., \cite{bertsekas1999nonlinear,bertsekas2003convex,nocedal2006numerical,boyd2004convex}) omit discussion of any of the related 
classes of operators that form the basis of our discussion; namely operators that are averaged, inverse strongly monotone, or pseudocontractive. 
We now briefly consider those references that do mention one (or, in some cases, two) of these operator classes.

In \cite{eckstein1989splitting} the discussion covers (projected) gradient descent, the proximal-point method, forward-backward splitting, 
the alternating direction method of multipliers, and numerous other methods. The coverage does touch on the special case of $1$-pseudocontractivity  (i.e., firm nonexpansiveness) and on the special case of $1$-inverse strong monotonicity; however, without the more general concepts of $\nu$-pseudocontractivity  and $\sigma$-inverse strong monotonicity (and the relationship between these concepts) there are significant limits to the results that can be obtained. One particularly notable limitation of a focus on firmly nonexpansive operators is that the class of firmly nonexpansive operators is not closed under composition.

\cite{vasin1995illposed} discusses both averaged operators and pseudocontractive operators, and also establishes the relationship between these classes; however, the  general concept of inverse strong monotonicity is not present; in addition to the absence of the combined relationship between pseudocontractivity and inverse strong monotonicity, there is no unified convergence rate analysis of the type we discuss. On the other hand, \cite{bauschke2011convex} discusses averagedness and inverse strong monotonicity (called there co-coercivity). Both averagedness and inverse strong monotonicity are discussed extensively and the latter is the central property on which most of the subsequent analysis is based. However, pseudocontractivity is essentially absent and convergence rates are not discussed.

\cite{golshtein1996modified} discusses inverse strong monotonicity extensively in their discussion of gradient-type methods and in methods involving modified monotone mappings; 
however, the concepts of averagedness and pseudocontractivity are absent. The discussion in \cite{byrne2004unified} establishes that many iterative optimization methods have updates that are \emph{averaged}; this work also includes the 
relationship between averaged and inverse strongly monotone. However, the concept of pseudocontractivity is absent. While the Krasnoselskii-Mann 
Theorem \cite{mann1953mean} is used to establish convergence of methods with averaged updates, there is no analysis of convergence rates.

Of more recent papers discussing proximal methods, such as \cite{combettes2005signal,combettes2011proximal}, the special case of firm nonexpansiveness appears, but not 
the more general cases of interest to us. The story is similar for other recent papers discussing the alternating direction method of multipliers, such as 
\cite{boyd2011distributed,parikh2013proximal,he2012rate}. This, again, limits the utility and breadth of applicability of the results obtained. 

In the machine learning community, as in the optimization community, the classes of operators that we consider pass largely unremarked. 
\cite{lee2012proximal} uses the proximal-point mapping and establishes the firm nonexpansiveness of some of the operations involved in their algorithm. The discussion in \cite{duchi2009efficient,singer2009efficient} squarely covers some of the methods that we touch on here; despite an apparent similarity in general outline, the analysis is pursued from a rather different perspective that is largely disjoint from our approach; in part this is due to a focus there on objective function value, 
in contrast to our focus on the sequence of iterates. One of the methods we discuss, the alternating direction method of multipliers, has recently been considered in stochastic and  online variants \cite{ouyang2013stochastic, wang2012online}, but our discussion here is restricted to the batch deterministic case.

\subsection{Preliminaries}
We restrict our attention to the finite dimensional case; this avoids situations in which the methods we discuss might attain weak convergence 
but not strong convergence \cite{guler1991convergence} (since in the finite-dimensional case, the concepts of weak convergence and strong convergence coincide).

We use the notation $\text{Fix }T$ to denote the fixed point set of the operator $T$; that is, $\text{Fix }T=\{x\in\mathbb{R}^n\ |\ T(x)=x\}$; we restrict our attention
throughout to operators having nonempty fixed point sets. We use $\left\Vert \cdot \right\Vert$ to denote the standard 2-norm and $\left\langle \cdot,\cdot\right\rangle$ to
denote the associated standard inner product. When convenient, we switch freely between operator notation and relation notation; for example, we might state 
the operator inclusion $y\in T(x)$ equivalently in relation form as $(x,y)\in T$.

With regard to terminology, we note that the property of inverse strong monotonicity \cite{rockafellar2004variational,golshtein1996modified} is elsewhere 
called co-coercivity \cite{combettes2009fejer}.
On the other hand, the term pseudocontractive is sometimes used to mean a significantly different property (cf \cite{bertsekas1989parallel}) than our usage here.

\section{Classes of Operators}

We begin with a number of definitions.

\begin{defn}
An operator   $ T : \mathbb{R}^n \to \mathbb{R}^n $ is \emph{nonexpansive} when
\begin{align}
\left\Vert Tx-Ty\right\Vert &\leq\left\Vert x-y\right\Vert 
 \end{align}
for all $x,y \in \mathbb{R}^n$.
\end{defn}

\begin{defn}
For any operator $ T : \mathbb{R}^n \to \mathbb{R}^n $, the associated \emph{displacement operator} 
is $G=I-T$.
\end{defn}

\begin{defn}
An operator   $ T : \mathbb{R}^n \to \mathbb{R}^n $ is \emph{firmly nonexpansive} when
\begin{align}
\left\Vert Tx-Ty\right\Vert ^{2}&\leq\left\Vert x-y\right\Vert ^{2}-\left\Vert Gx-Gy\right\Vert ^{2},
\end{align}
for all $x,y \in \mathbb{R}^n$, where $G$ is the displacement operator associated with $T$.
\end{defn}

\begin{defn}
An operator   $ T : \mathbb{R}^n \to \mathbb{R}^n $ is $\nu$-\emph{pseudocontractive} when
\begin{align}
\left\Vert Tx-Ty\right\Vert ^{2}&\leq\left\Vert x-y\right\Vert ^{2}-\nu\left\Vert Gx-Gy\right\Vert ^{2}
\end{align}
for all $x,y \in \mathbb{R}^n$, where $\nu >0$ and $G$ is the displacement operator associated with $T$.
\end{defn}

Note that firmly nonexpansive operators correspond to $\nu$-pseudocontractive operators in the special case of $\nu=1$.

\begin{defn}
An operator   $ T : \mathbb{R}^n \to \mathbb{R}^n $ is $\alpha$-\emph{averaged} when it is of the form
\begin{align}
T=(1-\alpha)I+\alpha N ,
\end{align}
where $N$ is a nonexpansive operator and $\alpha\in(0,1)$.
\end{defn}

We remark in passing that the use of the term averaged as defined above corresponds to the generalization 
from the special case in which the convex combination parameter $\alpha=\sfrac{1}{2}$.

\begin{defn}
A point-to-set operator   $T : \mathbb{R}^n \to 2^{\mathbb{R}^n} $ is \emph{monotone} when
\begin{align}
\left\langle u-v,x-y\right\rangle &\geq0
\end{align}
for any $(x,u)\in T$, $(y,v)\in T$, where $x,y\in\mathbb{R}^{n}$.
\end{defn}

\begin{defn}
An operator   $ T : \mathbb{R}^n \to \mathbb{R}^n $ is $\sigma$-\emph{inverse strongly monotone} when
\begin{align}
\left\langle Tx-Ty,x-y\right\rangle &\geq\sigma\left\Vert Tx-Ty\right\Vert ^{2}
\end{align}
for all $x,y\in\mathbb{R}^{n}$, where $\sigma>0$.
\end{defn}

\subsection{Pseudocontractive Updates: Prevalence and Common Convergence Rate Analysis}
We discuss the details in later sections after necessary theoretical developments; for now we 
state the two main themes of this manuscript.
\begin{prop*}[Details in Section \ref{sect:methods}]
The following prominent convex optimization methods have pseudocontractive updates:
\begin{enumerate}[label=\emph{\alph*})]
\item Orthogonal projection onto a nonempty closed convex set is $1$-pseudocontractive.
\item Over- or -under relaxed orthogonal projection methods are $(\frac{2-\omega}{\omega})$-pseudocontractive, where $\omega\in(0,2)$ is the relaxation parameter. 
\item For a convex function $f$ with an $L$-Lipschitz gradient, the gradient descent iterative 
update operator $I-\gamma\nabla f$ is $(\frac{2}{\gamma L}-1)$-pseudocontractive for any step factor $\gamma\in (0,\frac{2}{L})$. 
\item The proximal-point mapping, $\text{prox}_{f,\lambda}$,  associated with a convex function $f$ is $1$-pseudocontractive for any $\lambda>0$.
\item The $\omega$-relaxed proximal-point mapping, $\text{prox}_\omega$,  where $\omega\in(0,2)$, associated with 
a convex function $f$ is $(\frac{2-\omega}{\omega})$-pseudocontractive for any $\lambda>0$.
\item Projected gradient descent involving a convex function $f$ with an $L$-Lipschitz gradient is pseudocontractive 
with parameter $(1-\frac{\gamma L}{2})$, where $\gamma\in(0,\frac{2}{L})$ is the gradient descent step factor.
\item Over- or -under-relaxed projected gradient descent is pseudocontractive.
\item Forward-backward splitting in which the smooth convex function $f$ has an $L$-Lipschitz gradient is pseudocontractive 
with parameter $(1-\frac{\gamma L}{2})$, where $\gamma\in(0,\frac{2}{L})$ is the gradient descent step factor.
\item The most basic version of the alternating direction method of multipliers update is $1$-pseudocontractive. 
\end{enumerate}
\end{prop*}

In addition to highlighting the common role played by pseudocontractive updates in all of the methods above, 
we also present a novel proof establishing that the error criterion $\left\Vert x^{N}-Tx^{N}\right\Vert ^{2}$ is $o(1/N)$ for any method with a pseudocontractive update.

\begin{thm*}[Details in Section \ref{sect:proofs}]
The error criterion $\left\Vert x^{N}-Tx^{N}\right\Vert ^{2}$ is $o(1/N)$ for any method with a pseudocontractive iteration operator. 
\end{thm*}

Having stated our primary results, we now begin introducing the background concepts that we will use in establishing these results.

\subsection{Relationships}
We note that the following statements are equivalent: 
\begin{itemize}
 \item $T$ is  $\left(\frac{1-\alpha}{\alpha}\right)$-pseudocontractive   
 
 \item $T$ is $\alpha$-averaged
 
 \item $I-T$ is $\frac{1}{2\alpha }$-inverse strongly monotone
 
\end{itemize}

The relationship between pseudocontractive and averaged can be found in \cite{vasin1995illposed}. 
The relationship between averaged and inverse strongly monotone can be found in \cite{byrne2004unified}. 
The combination of all three of the concepts of pseudocontractivity, averagedness, and inverse strong monotonicity in the present manuscript 
provides the foundation for most of our results.

We also note an equivalence in a special case: $T$ being $1$-pseudocontractive is equivalent to $T$ being $1$-inverse strongly monotone;
see Corollary \ref{cor:equivalence-fne-and-ism}. 

The significance of these relationships is that they allow us to first establish properties in whatever form is most immediate
and then to translate the properties to the most convenient form for subsequent use.

\subsubsection{Over- or under-relaxation involving a firmly nonexpansive operator}
Finally, we note that if we start with a firmly nonexpansive operator $F$, it is actually the case that the operator $T=(1-\omega)I+\omega F$ is averaged
not only for $\omega\in(0,1)$, but for $\omega\in(0,2)$. This is because when $F$ is firmly nonexpansive, its 
associated reflection\footnote{The use of the term reflection for this operator is most immediately clear when we consider
the case of projection onto a hyperplane $H$. Denoting projection onto the hyperplane $H$ by $P_H$ and denoting reflection
across the hyperplane by $R_H$, we observe that $P_H=\sfrac{1}{2}I+\sfrac{1}{2}R_H$. The usage in the case above is a natural generalization.} 
operator $R=2F-I$ is nonexpansive
(see Proposition \ref{prop:reflection-nonexpansive}). 

We state this result on relaxation in the following Proposition for future reference:
\begin{prop}
\label{prop:relaxation-with-fne}
For a firmly nonexpansive operator $F$, the operator $T=(1-\omega)I+\omega F$ is $(\frac{2-\omega}{\omega})$-pseudocontractive for $\omega\in(0,2)$.
\end{prop}
\begin{prf}
Denoting the reflection operator associated with $F$ as $R=2F-I$, we note that 
\begin{align}
\label{eqn:relaxation}
\left(1-\omega\right)I+\omega F&=\left(1-\frac{\omega}{2}\right)I+\frac{\omega}{2}R. 
\end{align}
Since Proposition \ref{prop:reflection-nonexpansive} establishes that $R$ is nonexpansive,
the expression above means that the operator  $T=(1-\omega)I+\omega F$ can be seen to be
$\frac{\omega}{2}$-averaged for $\omega\in(0,2)$; this in turn corresponds to 
$T$ being $(\frac{2-\omega}{\omega})$-pseudocontractive for $\omega\in(0,2)\square$. 
\end{prf}

The significance of the result above is that iterations involving under- or over-relaxation of a pseudocontractive update operator are also covered by our results.
Over-relaxation, in which $\omega\in(1,2)$, is a broadly useful technique most commonly used to accelerate convergence of iterative methods in numerical linear
algebra \cite{varga2009matrix}; the observation that over-relaxation of a firmly nonexpansive operator yields a pseudocontractive update is thus of some notable interest.

\subsection{Closure Properties}
The convex combination of nonexpansive operators is a nonexpansive operator.
The composition of nonexpansive operators is a nonexpansive operator.
The convex combination of pseudocontractive operators is a pseudocontractive operator.
The composition of pseudocontractive operators is a pseudocontractive operator.
These results can be found in \cite{vasin1995illposed}; however, we can in fact be even more specific: 
\begin{thm}[\cite{ogura2002nonstrictly} Theorem 3]
\label{thm:closure-results-pseudocontractive}
When $T_{1}$ is $\alpha_{1}$-averaged, where $\alpha_{1}\in\left[0,1\right)$, and $T_{2}$ is $\alpha_{2}$-averaged, where $\alpha_{2}\in\left[0,1\right)$, we have
\begin{enumerate}[label=\emph{\alph*})]
\item for $t\in[0,1]$, the $t$-convex combination $T_{\text{cc}}=(1-t) T_{1}+tT_{2}$ is $\alpha_{\text{cc}}$-averaged, 
where $\alpha_{\text{cc}}=(1-t)\alpha_1+t\alpha_2$, 
\item the composition $T_{\text{co}}=T_{1}T_{2}$ is $\alpha_{\text{co}}$-averaged, 
where  $\alpha_{\text{co}}=\frac{\alpha_{1}+\alpha_{2}-2\alpha_{1}\alpha_{2}}{1-\alpha_{1}\alpha_{2}}$.
\end{enumerate}
\end{thm}

Using the relationship between averaged and pseudocontractive, we note that the convex combination $T_\text{cc}$ is $\nu_{\text{cc}}$-pseudocontractive, where $\nu_{\text{cc}}=\frac{1-\alpha_{\text{cc}}}{\alpha_{\text{cc}}}$; similarly,
the composition $T_{\text{co}}$ is $\nu_{\text{co}}$-pseudocontractive, 
where $\nu_{\text{co}}=\frac{1-\alpha_{1}-\alpha_{2}+\alpha_{1}\alpha_{2}}{\alpha_{1}+\alpha_{2}-2\alpha_{1}\alpha_{2}}.$

Note that while this result establishes that both $T_{1}T_{2}$ and $T_{2}T_{1}$ are pseudocontractive, 
the respective fixed point sets need not coincide\footnote{For example, consider the case of successive orthogonal projection onto two closed, convex, nonempty sets 
$\mathcal{A}$ and $\mathcal{B}$, with $\mathcal{A}\cap\mathcal{B}=\emptyset$. In this situation, 
$P_{\mathcal{A}}P_{\mathcal{B}}$ and $P_{\mathcal{B}}P_{\mathcal{A}}$ each have nonempty fixed point sets, but $\text{Fix }P_{\mathcal{A}}P_{\mathcal{B}}$ is 
the set of points in $\mathcal{A}$ that are at the minimal distance to $\mathcal{B}$, whereas $\text{Fix }P_{\mathcal{B}}P_{\mathcal{A}}$ is the set of 
points in $\mathcal{B}$ that are at the minimal distance to $\mathcal{A}$.}.

The significance of these closure properties is twofold: first, closure under composition and convex combination means that new algorithms can very naturally
be constructed by using previous algorithms as building blocks; second, the resulting constructed algorithms will maintain the desirable properties that 
the individual building block updates possessed (and using Theorem \ref{thm:closure-results-pseudocontractive} we can track the associated pseudocontractivity
parameter values from the simple methods to the more complex combined methods).

\section{Methods with Pseudocontractive Updates}
\label{sect:methods}
In this section we demonstrate that pseudocontractive updates are widespread in many popular convex optimization methods.

\subsection{Projection}
Most discussions of projection methods only establish that projection is nonexpansive; the following known result 
demonstrates that projection is in fact $1$-pseudocontractive (i.e., firmly nonexpansive).

\begin{prop}
Orthogonal projection $P_\mathcal{C}$ onto a nonempty closed convex set $\mathcal{C}\subset\mathbb{R}^n$ is $1$-pseudocontractive. For convenience we show the equivalent 
result that projection is $1$-inverse strongly monotone; that is, 
\begin{align}
\left\langle P_{\mathcal{C}}x-P_{\mathcal{C}}y,x-y\right\rangle &\geq\left\Vert P_{\mathcal{C}}x-P_{\mathcal{C}}y\right\Vert ^{2}, 
\end{align}
for all $x,y\in\mathbb{R}^n$.

\end{prop}
\begin{prf}
The variational characterization of projection is: $z$ is the projection of $x$ onto $\mathcal{C}$ if and only if $\left\langle c-z,x-z\right\rangle \leq0$ for every $c\in\mathcal{C}$
(see, e.g., \cite{ruszczynski2006nonlinear}); we will express this as 
\begin{align}
\left\langle c-P_{\mathcal{C}}x,x-P_{\mathcal{C}}x\right\rangle &\leq0, 
\end{align}
for every $c\in\mathcal{C}$.
From the preceding characterization, we have, for any $x,y\in \mathbb{R}^n$ 
\begin{align}
\left\langle P_{\mathcal{C}}x-P_{\mathcal{C}}y,y-P_{\mathcal{C}}y\right\rangle &\leq0 \\
\left\langle P_{\mathcal{C}}y-P_{\mathcal{C}}x,x-P_{\mathcal{C}}x\right\rangle &\leq0,
\end{align}
yielding the immediate conclusion
\begin{align}
\left\Vert P_{\mathcal{C}}x-P_{\mathcal{C}}y\right\Vert ^{2}&\leq\left\langle P_{\mathcal{C}}x-P_{\mathcal{C}}y,x-y\right\rangle, 
\end{align}
for any $x,y\in\mathbb{R}^n$.
The result immediately above establishes that $P_{\mathcal{C}}$ is $1$-inverse strongly monotone, so that Corollary \ref{cor:equivalence-fne-and-ism} implies 
that  $P_{\mathcal{C}}$ is $1$-pseudocontractive. $\square$
\end{prf}

\subsection{Over- or Under-relaxed Projection}
Using the notation $P$ for orthogonal projection onto the nonempty closed convex set $\mathcal{C}$, we introduce the $\omega$-relaxed projection operator 
$P_{\omega}=(1-\omega)I+\omega P$, where $\omega\in(0,2)$. Since we just established that $P$ is $1$-pseudocontractive, we conclude that the 
associated reflection operator $2P-I$ is nonexpansive, and so from our previous reasoning in Equation \ref{eqn:relaxation} we observe that
$P_{\omega}$ is $(\frac{2-\omega}{\omega})$-pseudocontractive.

We note in passing that reasoning similar to the reasoning in the previous two subsections can also be used to establish 
that cyclic subgradient projection (CSP) methods (in both their relaxed and unrelaxed forms) \cite{censor1997parallel} involve pseudocontractive updates.

\subsection{Gradient Descent}
Typical discussions of gradient descent do not make any reference to operator theory, or, more specifically, to pseudocontractivity; the known result below 
establishes conditions under which the gradient descent iterative update is pseudocontractive.
The explicit form of the gradient descent iterative update is 
\begin{align}
x^{k+1}&=\left[I-\gamma\nabla f\right]\left(x^{k}\right).
\end{align}
\begin{prop}
For a convex function $f$ with an $L$-Lipschitz gradient, the gradient descent iterative 
update operator $I-\gamma\nabla f$ is $(\frac{2}{\gamma L}-1)$-pseudocontractive for any step factor $\gamma\in (0,\frac{2}{L})$. 
\end{prop}
\begin{prf}
We have assumed that $\nabla f$ is $L$-Lipschitz, which implies that $\frac{1}{L}\nabla f$ is nonexpansive.  
Theorem 6.9 from Chapter 1 of \cite{golshtein1996modified} establishes that $\frac{1}{L}\nabla f$ nonexpansive implies that $\frac{1}{L}\nabla f$ is $1$-pseudocontractive.

Corollary \ref{cor:equivalence-fne-and-ism} establishes that $\frac{1}{L}\nabla f$ being $1$-pseudocontractive is 
equivalent to $\frac{1}{L}\nabla f$ being $1$-inverse strongly monotone. 
Lemma \ref{lem:ism-scaling} establishes that $\frac{1}{L}\nabla f$ being $1$-inverse strongly monotone implies that 
 $\gamma\nabla f$ is $\frac{1}{\gamma L}$-inverse strongly monotone, since $\gamma L>0$.

When an operator $T$  is exactly $\sfrac{1}{2}$-inverse strongly monotone, $G=I-T$ is nonexpansive.
When $T$ is $\sigma$-inverse strongly monotone for $\sigma\in\left(\sfrac{1}{2},+\infty\right)$, $G=I-T$ is $(2\sigma-1)$-pseudocontractive.
To ensure $I-\gamma\nabla f$ is pseudocontractive, we must thus ensure that $\gamma\nabla f$  is $\sigma$-inverse strongly monotone 
with $\sigma\in\left(\sfrac{1}{2},+\infty\right)$; equivalently, we must ensure that $(2\sigma-1)>0$.

$\gamma\nabla f$ being $\frac{1}{\gamma L}$-inverse strongly monotone corresponds to $I-\gamma\nabla f$
being $\left(\frac{2}{\gamma L}-1\right)$-pseudocontractive. 
The pseudocontractivity parameter is required to be strictly positive; in order to ensure that
$\left(\frac{2}{\gamma L}-1\right)>0$,  we require $\gamma\in\left(0,\frac{2}{L}\right)$,
since $\left(\frac{2}{\gamma L}-1\right)>0$ requires $\frac{2}{L}>\gamma$. $\square$
\end{prf}

\subsection{Proximal-Point Method}
As is the case with projection, most discussions of proximal-point methods only establish that the proximal-point mapping is nonexpansive; the following
known result demonstrates that the proximal-point mapping is in fact $1$-pseudocontractive (i.e., firmly nonexpansive).
We note that the proximal-point mapping associated with $f$ can be written in either of the two equivalent forms
\begin{align}
\text{prox}_{f,\lambda}(x)&=\underset{z\in\mathbb{R}^{n}}{\text{argmin }}\left\{ f\left(z\right)+\frac{1}{2\lambda}\left\Vert z-x\right\Vert ^{2}\right\} \\
\text{prox}_{f,\lambda}(x)&=(I+\lambda\partial f)^{-1}(x),
\end{align}
where $\lambda>0$; see, e.g., \cite{rockafellar1970convex}. The mapping $(I+\lambda\partial f)^{-1}$ is called the $\lambda$-resolvent of the subdifferential mapping $\partial f$.

\begin{prop}
\label{prop:prox-fne}
For a convex function $f$, the associated proximal-point update is $1$-pseudocontractive, for any $\lambda>0$.
\end{prop}
\begin{prf}
Consider  $x^{+}=\text{prox}_{f,\lambda}(x)$ and $y^{+}=\text{prox}_{f,\lambda}(y)$. 
Lemma \ref{lem:prox-inclusion} tells us that 
\begin{align}
\left(x^{+},x-x^{+}\right)&\in\lambda\partial f\\
\left(y^{+},y-y^{+}\right)&\in\lambda\partial f.
\end{align}
The subdifferential mapping of a convex function is monotone \cite{rockafellar1970convex}; thus $\lambda\partial f$, where $\lambda>0$, is also monotone. From this we have
\begin{align}
\left\langle x^{+}-y^{+},\left( x-x^{+}\right) -\left( y-y^{+}\right) \right\rangle &\geq0 \\
\left\langle x^{+}-y^{+},\left( x-y\right) -\left( x^{+}-y^{+}\right) \right\rangle &\geq0 \\
\left\langle x^{+}-y^{+},x-y\right\rangle &\geq\left\Vert x^{+}-y^{+}\right\Vert ^{2}.
\end{align}
This establishes that $\text{prox}_{f,\lambda}$  is $1$-inverse strongly monotone; applying Corollary \ref{cor:equivalence-fne-and-ism} then establishes that $\text{prox}_{f,\lambda}$ 
is $1$-pseudocontractive for any $\lambda>0$. $\square$
\end{prf}

Although we do not go into further detail, we note here that analogous arguments apply to any algorithm 
involving the resolvent operator associated with a maximal monotone mapping (cf \cite{eckstein1989splitting}).

\subsection{Relaxed proximal-point Method}
Lightening our notation and for the moment taking $\text{prox}$ to denote the $\lambda$-proximal-point mapping associated with the convex function $f$, 
we introduce the $\omega$-relaxed proximal-point iteration operator
$\text{prox}_{\omega}=(1-\omega)I+\omega \text{prox}$, where $\omega\in(0,2)$. Since we have established that $\text{prox}$ is $1$-pseudocontractive, we conclude that the 
associated reflection operator $2\text{prox}-I$ is nonexpansive, and so from our previous reasoning in Equation \ref{eqn:relaxation} we observe that
$\text{prox}_{\omega}$ is $(\frac{2-\omega}{\omega})$-pseudocontractive.

\subsection{Projected Gradient Descent}
The iterative update for projected gradient descent involving a convex function $f$ with an $L$-Lipschitz gradient and a nonempty closed convex constraint set $\mathcal{C}$
corresponds to the composition $P_{\mathcal{C}}\circ[I-\gamma \nabla f]$. 
We have established that the projection $P_{\mathcal{C}}$ is $1$-pseudocontractive and that the gradient descent update $I-\gamma \nabla f$ 
is $(\frac{2}{\gamma L}-1)$-pseudocontractive, where $\gamma\in(0,\frac{2}{L})$. 
Since the composition of pseudocontractive updates is pseudocontractive \cite{vasin1995illposed}, we have thus established that projected gradient descent 
is pseudocontractive. We also note that under- or over-relaxation of the projection portion of the update still yields
pseudocontractivity. 

We can use Theorem \ref{thm:closure-results-pseudocontractive}b to be more specific.  
Basic projected gradient descent is the composition of the gradient descent step, which is $\frac{\gamma L}{2}$-averaged, and the projection step, which is $\sfrac{1}{2}$-averaged; 
the composition thus is $(1-\frac{\gamma L}{2})$-pseudocontractive.
Under- or over-relaxed projected gradient descent has the $\sfrac{1}{2}$-averaged projection replaced 
with $\omega$-relaxed projection, which we have seen is $(\frac{2-\omega}{\omega})$-pseudocontractive.

\subsection{Forward-Backward Splitting}
Consider the optimization problem 
\begin{align}
\underset{x\in\mathbb{R}^{n}}{\text{minimize }}&f\left(x\right)+g\left(x\right),
\end{align}
where $f$ and $g$ are each convex and $f$ has an $L$-Lipschitz gradient.
The explicit statement of the forward-backward splitting method is 
\begin{align}
x^{k+1}=\text{prox}_{g,\lambda}\circ[I-\gamma \nabla f](x^k) 
\end{align}

When used as a gradient descent step followed by a proximal step, so long as the gradient descent step factor is $\gamma\in (0,\frac{2}{L})$, 
we have the composition of pseudocontractive operators, which is in turn pseudocontractive.

We can again use Theorem \ref{thm:closure-results-pseudocontractive}b to be more specific.  
Basic forward-backward splitting is the composition of the gradient descent step, which is $\frac{\gamma L}{2}$-averaged, and the proximal step, which is $\sfrac{1}{2}$-averaged; 
the composition thus is $(1-\frac{\gamma L}{2})$-pseudocontractive.
When we use under- or over-relaxation of the proximal step in forward-backward splitting, the $\sfrac{1}{2}$-averaged proximal step is replaced 
with an $\omega$-relaxed generalized proximal step, which we have seen is $(\frac{2-\omega}{\omega})$-pseudocontractive.

\subsection{Alternating Direction Method of Multipliers}
Consider the optimization problem 
\begin{align}
\underset{x\in\mathbb{R}^{n}}{\text{minimize }}&f\left(x\right)+g\left(x\right),
\end{align}
where $f$ and $g$ are each convex.
Using the notation $R_f=2\text{prox}_{f,\lambda}-I$ and $R_g=2\text{prox}_{g,\lambda}-I$,
we can express \cite{boyd2011splitting} the most basic version of the alternating direction method of multipliers (ADMM) update in terms of operators as
\begin{align}
x^{k+1}=[\sfrac{1}{2}I+\sfrac{1}{2} R_g R_f](x^k). 
\end{align}

Proposition \ref{prop:prox-fne} establishes that $\text{prox}_{f,\lambda}$ and $\text{prox}_{g,\lambda}$ are each $1$-pseudocontractive, for any $\lambda>0$.
Proposition \ref{prop:reflection-nonexpansive} then establishes that the respective reflection operators $R_f$ and $R_g$ are each nonexpansive.
The composition of nonexpansive operators is nonexpansive, so $R_g R_f$ is nonexpansive. We thus observe that the most basic version of the ADMM update 
is $\sfrac{1}{2}$-averaged, and thus is $1$-pseudocontractive. We also immediately observe that the more general convex-combination ADMM-like update
\begin{align}
x^{k+1}=[(1-\alpha) I+\alpha R_g R_f](x^k), 
\end{align}
where $\alpha\in(0,1)$, will be $(\frac{1-\alpha}{\alpha})$-pseudocontractive.

\section{Convergence Rate of Methods with $\nu$-Pseudocontractive Updates}
\label{sect:proofs}
The worst-case convergence rate is well-known for each of the methods that we discuss in Section \ref{sect:methods}. 
The standard approach generally requires separate arguments to establish this rate for each method; in contrast, by leveraging the fact that
the updates in each of these methods are pseudocontractive (and making use of the relationships between pseudocontractive, averaged, and inverse strongly monotone)
we are able to use one argument for all of these methods.

Consider a method for which the iteration operator $T$ is $\nu$-pseudocontractive. We will establish 

\begin{lem}
\label{lem:monotone-decreasing}
For a $\nu$-pseudocontractive iteration operator $T$, the norm of the displacement operator $G=I-T$ satisfies the expression 
\begin{align}
\left\Vert Gx^{k+1}\right\Vert ^{2}&\leq\left\Vert Gx^{k}\right\Vert ^{2}-\nu\left\Vert Gx^{k}-Gx^{k+1}\right\Vert ^{2} 
\end{align} 
for $k\in\{0,1,\ldots\}$.
\end{lem}
\begin{prf}
From the relationship between pseudocontractivity and inverse strong monotonicity, we note that when $T$ is $\nu$-pseudocontractive, 
we have that $G=I-T$ is $\frac{1+\nu}{2}$-inverse strongly monotone. That is, we have
\begin{align}
\left\langle Gx-Gy,x-y\right\rangle &\geq\left(\frac{1+\nu}{2}\right)\left\Vert Gx-Gy\right\Vert ^{2},
\end{align}
for all $x,y\in \mathbb{R}^n$.
Thus,
\begin{align}
\left\langle Gx^{k}-Gx^{k+1},x^{k}-x^{k+1}\right\rangle &\geq\frac{1+\nu}{2}\left\Vert Gx^{k}-Gx^{k+1}\right\Vert ^{2} \\
\left\Vert Gx^{k}\right\Vert ^{2}-\left\langle Gx^{k+1},Gx^{k}\right\rangle &\geq\frac{1+\nu}{2}\left\Vert Gx^{k}-Gx^{k+1}\right\Vert ^{2},
\end{align}
which yields
\begin{align}
\left\Vert Gx^{k+1}\right\Vert ^{2}&\leq\left\Vert Gx^{k}\right\Vert ^{2}-\nu\left\Vert Gx^{k}-Gx^{k+1}\right\Vert ^{2},
\end{align}
where we have used $Gx^k=x^{k}-x^{k+1}$ and
$2\left\langle Gx^{k+1},Gx^{k}\right\rangle =\left\Vert Gx^{k+1}\right\Vert ^{2}+\left\Vert Gx^{k}\right\Vert ^{2}-\left\Vert Gx^{k+1}-Gx^{k}\right\Vert ^{2}$. $\square$
\end{prf}

Using Lemma \ref{lem:monotone-decreasing}, we can now show that 

\begin{thm} 
\label{thm:convergence-rate}
Any method with a $\nu$-pseudocontractive iteration operator has error criterion satisfying $\left\Vert x^{N}-Tx^{N}\right\Vert ^{2}\leq\frac{1}{\left(N+1\right)\nu}\left\Vert x^{0}-x^{\star}\right\Vert ^{2} $ at iteration $N$. 
\end{thm}
\begin{prf}
An alternate expression of $\nu$-pseudocontractivity with one argument, $x^\star$, coming from the fixed point set is 
\begin{align}
 \nu\left\Vert Gx^{k}\right\Vert^{2}&\leq\left\Vert x^{k}-x^{\star}\right\Vert^{2}-\left\Vert x^{k+1}-x^{\star}\right\Vert^{2}.
\end{align}
Summing over all the iterations from $0$ to $N$ yields
\begin{align}
\sum_{k=0}^{N}\left\Vert Gx^{k}\right\Vert ^{2}&\leq\frac{1}{\nu}\left\Vert x^{0}-x^{\star}\right\Vert ^{2}.
\end{align}
Lemma \ref{lem:monotone-decreasing} tells us that $\left\Vert Gx^{k+1}\right\Vert^{2} \leq\left\Vert Gx^{k}\right\Vert^{2}$ for each $k$, so we conclude that 
\begin{align}
\left\Vert Gx^{N}\right\Vert ^{2}\leq&\frac{1}{\left(N+1\right)\nu}\left\Vert x^{0}-x^{\star}\right\Vert ^{2}.\  \square
\end{align}

\end{prf}

We can establish an even stronger result by leveraging a useful technical Lemma from \cite{dong2013proximal}.

\begin{lem} [\cite{dong2013proximal}]
\label{lem:sequence-lemma}
Consider three sequences of strictly positive numbers $\{\alpha_k\},\{\beta_k\},\{\gamma_k\}$. When these sequences are such that $\{\beta_k\}$
is unsummable, $\{\gamma_k\}$ is nonincreasing, and the relationship $\alpha_{k+1}^{2}\leq\alpha_{k}^{2}-\beta_{k}\gamma_{k}$ holds for each $k\in\left\{ 0,1,\ldots\right\}$, then it will also be the case that there exists another sequence $\{\varepsilon_k\}$ such that 
\begin{itemize}
 \item $\varepsilon_{k}\in\left(\alpha_{k},\alpha_{0}\right),$ 
 \item $\lim_{k\to+\infty}\varepsilon_{k}=\lim_{k\to+\infty}\alpha_{k},$
 \item $\gamma_{k}\cdot\left[\sum_{i=0}^{k}\beta_{i}\right]\leq2\cdot\alpha_{0}\cdot\varepsilon_{k}.$
\end{itemize}

\end{lem}

We now use this Lemma to establish
\begin{thm} 
\label{thm:convergence-rate-stronger}
For any method with a $\nu$-pseudocontractive iteration operator, the error criterion $\left\Vert x^{N}-Tx^{N}\right\Vert ^{2}$ is $o(1/N)$. 
\end{thm}
\begin{prf}
We will apply Lemma \ref{lem:sequence-lemma} with $\alpha_k=\left\Vert x^{k}-x^{\star}\right\Vert$, $\beta_k=\nu$, and $\gamma_k=\left\Vert Gx^{k}\right\Vert ^{2}$. An infinite sequence of constants is clearly unsummable and Lemma \ref{lem:monotone-decreasing} establishes that $\{\left\Vert Gx^{k}\right\Vert ^{2}\}$ is nonincreasing. 

Moreover, the definition of $\nu$-pseudocontractivity when one argument, $x^{\star}$, comes from the fixed point set provides us with the necessary relationship 
\begin{align}
\left\Vert x^{k+1}-x^{\star}\right\Vert ^{2}\leq\left\Vert x^{k}-x^{\star}\right\Vert ^{2}-\nu\left\Vert Gx^{k}\right\Vert ^{2}  
\end{align}
for each $k\in\left\{ 0,1,\ldots\right\}.$ 

Thus, Lemma \ref{lem:sequence-lemma} tells us that there exists a sequence  $\{\varepsilon_k\}$ such that
\begin{itemize}
 \item $\varepsilon_{k}\in\left(\left\Vert x^{k}-x^{\star}\right\Vert ,\left\Vert x^{0}-x^{\star}\right\Vert \right),$ 
 \item $\lim_{k\to+\infty}\varepsilon_{k}=\lim_{k\to+\infty}\left\Vert x^{k}-x^{\star}\right\Vert,$ 
 \item $\left\Vert Gx^{k}\right\Vert ^{2}\cdot\left[\sum_{i=0}^{k}\nu\right]\leq2\cdot\left\Vert x^{0}-x^{\star}\right\Vert \cdot\varepsilon_{k}.$
\end{itemize}

Taking limsup on both sides of the preceding inequality  and noting that $\lim_{k\to+\infty}\varepsilon_{k}$
  exists and is equal to $\lim_{k\to+\infty}\left\Vert x^{k}-x^{\star}\right\Vert$ 
  we have 
  \begin{align}
  \limsup_{k\to\infty}\left\{ \left\Vert Gx^{k}\right\Vert ^{2}\left[\sum_{i=0}^{k}\nu\right]\right\} &\leq\lim_{k\to\infty}\left\{ 2\left\Vert x^{0}-x^{\star}\right\Vert \varepsilon_{k}\right\}
  \end{align}
  \begin{align}
  \limsup_{k\to\infty}\left\{ \left\Vert Gx^{k}\right\Vert ^{2}\left[\nu\left(k+1\right)\right]\right\} &\leq \nonumber \\
  &2\left\Vert x^{0}-x^{\star}\right\Vert \lim_{k\to \infty}\left\Vert x^{k}-x^{\star}\right\Vert
  \end{align}
  \begin{align}
  \limsup_{k\to\infty}\left\{ \left\Vert Gx^{k}\right\Vert ^{2}\left[\nu\left(k+1\right)\right]\right\} &\leq0,
  \end{align}
where we have used $\lim_{k\to+\infty}\left\Vert x^{k}-x^{\star}\right\Vert=0$ (for our finite dimensional setting, this follows because, e.g., the Krasnoselskii-Mann 
Theorem \cite{mann1953mean} establishes weak convergence of methods with averaged updates [and thus with pseudocontractive updates] and weak convergence coincides with strong convergence in finite dimension).

The final inequality above establishes the desired result: in terms of iteration $N$, we observe that  $\left\Vert Gx^{N}\right\Vert ^{2}$  is $o(1/N).\square$
\end{prf}

\section{Discussion}
Pseudocontractivity is a property largely overlooked by the machine learning community.  This paper takes a first step toward showing the benefits of recognizing and using pseudocontractivity. We have demonstrated the broad prevalence of pseudocontractive updates in existing optimization methods. We also emphasize that optimization methods made up of updates derived from compositions or convex combinations of pseudocontractive operators will also be pseudocontractive --- this provides a perspective on the convergence behavior of combined methods that appears to be entirely absent from classical analyses. We provide a demonstration of a unified analysis that applies to all pseudocontractive methods. 

With an eye to following these ideas from theory to practice, we are extending the analysis in this paper to the settings of inexact and stochastic pseudocontractive updates. Recent computational experience suggests that allowing inexact updates can have provide a significant boost to the practical efficiency of many optimization methods. We also anticipate that a Nesterov-style estimate sequence approach will provide a basis for accelerated versions of any pseudocontractive method. We believe that all of these characteristics provide strong motivation for researchers to examine optimization methods for unlooked-for pseudocontractivity.

\appendix
\section{Appendix}
The following known results are necessary for establishing the results elsewhere in the manuscript. 

\begin{prop}
\label{prop:fne-displacement}
An operator $F$ is $1$-pseudocontractive if and only if the displacement operator $G=I-F$ is $1$-pseudocontractive.
\end{prop}
\begin{prf}
Immediate from the definition of $1$-pseudocontractive. $\square$ 
\end{prf}

\begin{prop}
\label{prop:ism-displacement}
An operator $F$ is $1$-inverse strongly monotone if and only if the displacement operator $G=I-F$ is $1$-inverse strongly monotone.
\end{prop}
\begin{prf}
The relevant expressions for either direction are
\begin{align}
\left\Vert Gx-Gy\right\Vert ^{2}&\leq\left\langle Gx-Gy,x-y\right\rangle \\
\left\Vert \left(x-y\right)-\left(Fx-Fy\right)\right\Vert ^{2}&\leq\left\langle x-y,x-y\right\rangle \\
&\quad-\left\langle Fx-Fy,x-y\right\rangle \notag\\
\left\Vert Fx-Fy\right\Vert ^{2}&\leq\left\langle x-y,Fx-Fy\right\rangle. \ \square 
\end{align}

\end{prf}

\begin{prop}
\label{prop:equivalence-fne-and-displacement-ism}
An operator $F$ is $1$-pseudocontractive if and only if the displacement mapping $G=I-F$ is $1$-inverse strongly monotone.
\end{prop}
\begin{prf}
The relevant expressions for either direction are
\begin{align}
\left\Vert Tx-Ty\right\Vert ^{2}&\leq\left\Vert x-y\right\Vert ^{2}-\left\Vert Gx-Gy\right\Vert ^{2} \\
\left\Vert Tx-Ty\right\Vert ^{2}&+\left\Vert \left[x-y\right]-\left[Tx-Ty\right]\right\Vert ^{2}\leq\left\Vert x-y\right\Vert ^{2} \\
\left\Vert Tx-Ty\right\Vert ^{2}&\leq\left\langle x-y,Tx-Ty\right\rangle.\ \square 
\end{align}
\end{prf}

\begin{cor}
\label{cor:equivalence-fne-and-ism}
An operator $F$ is $1$-pseudocontractive if and only if $F$ is $1$-inverse strongly monotone.
\end{cor}
\begin{prf}
Proposition \ref{prop:equivalence-fne-and-displacement-ism} shows that $F$ is $1$-pseudocontractive if and only if $G=I-F$ is $1$-inverse strongly monotone.
Proposition \ref{prop:ism-displacement} shows that $G=I-F$ is $1$-inverse strongly monotone if and only if $F$ is $1$-inverse strongly monotone. 
Together these establish the result. $\square$
\end{prf}

\begin{prop}[\cite{goebel1990topics}, Theorem 12.1]
\label{prop:reflection-nonexpansive}
For a firmly nonexpansive operator $F$, the reflection operator $R=2F-I$ is nonexpansive.
\end{prop}
\begin{prf}
\begin{align}
\left\Vert Rx-Ry\right\Vert ^{2}&=4\left\Vert Fx-Fy\right\Vert ^{2}+\left\Vert x-y\right\Vert ^{2}\\
&\quad -4\left\langle Fx-Fy,x-y\right\rangle \notag\\
&\leq\left\Vert x-y\right\Vert ^{2},
\end{align}
for all $x,y\in\mathbb{R}^n$, where we have used the fact that $F$ firmly nonexpansive implies that $F$ is $1$-inverse strongly monotone; that is, 
that $\left\langle Fx-Fy,x-y\right\rangle \geq\left\Vert Fx-Fy\right\Vert ^{2}$, for all $x,y\in\mathbb{R}^n$. $\square$
\end{prf}

\begin{lem}
\label{lem:prox-inclusion}
For a convex function $f$, the proximal-point $x^{+}=\text{prox}_{f,\lambda}(x)$ satisfies the inclusion
\begin{align}
x-x^{+}\in \lambda \partial f(x^{+}), 
\end{align}
where $\lambda>0$.
\end{lem}
\begin{prf}
The resolvent statement of the proximal-point mapping is 
\begin{align}
\text{prox}_{f,\lambda}(x)&=(I+\lambda\partial f)^{-1}(x).
\end{align}
From $x^{+}=\text{prox}_{f,\lambda}(x)$, we thus have
\begin{align}
\left(x,x^{+}\right)&\in(I+\lambda\partial f)^{-1} \\
\left(x^{+},x\right)&\in I+\lambda\partial f \\
x-x^{+}&\in\lambda\partial f\left(x^{+}\right),
\end{align}
where we have freely shifted our view from relation to operator as convenient. $\square$
\end{prf}

\begin{lem}
\label{lem:ism-scaling}
$T$ being $\sigma$-inverse strongly monotone implies that, for $\gamma>0$, the scaled operator $\gamma T$ is $\frac{\sigma}{\gamma}$-inverse strongly monotone.
\end{lem}
\begin{prf}
Immediate from the definition of inverse strongly monotone. $\square$ 
\end{prf}

\bibliography{short-comments}
\bibliographystyle{icml2014}

\end{document}